\documentclass[12pt,twoside]{amsart}
\usepackage{amssymb, amsmath, amscd}
\newtheorem{thm}{Theorem}[section]
\newtheorem{crl}[thm]{Corollary}%[section]
\newtheorem{lmm}[thm]{Lemma}%[section]
\newtheorem{prp}[thm]{Proposition}%[section]
\theoremstyle{definition}
\newtheorem{dfn}[thm]{Definition}%[section]
\newtheorem{exa}[thm]{Example}%%%
\newtheorem{rem}[thm]{Remark}%%%
\newcommand{\xdiscrep}[0]{{\operatorname{discrep}}}
\newcommand{\xSupp}[0]{{\operatorname{Supp}}}
\newcommand{\Exc}[0]{{\operatorname{Exc}}}
\newcommand{\codim}[0]{{\operatorname{codim}}}

\newcommand{\Center}[0]{{\operatorname{Center}}}
\newtheorem*{ack}{{\bf{Acknowledgments}}}         %\renewcommand{\theack}{} 
 
\title[Addendum]
{Addendum to "Termination of $4$-fold canonical flips"}
\author{Osamu Fujino}
\date{2004/8/7} 
\address{Graduate School of Mathematics, 
Nagoya University, Chikusa-ku Nagoya 464-8602, Japan}
\email{fujino@math.nagoya-u.ac.jp}
\begin{document}
\bibliographystyle{amsalpha+}

\begin{abstract}
The definition of the weighted version of difficulty in 
"Termination of $4$-fold canonical flips" 
contains some errors. 
In this paper, we describe these errors 
and how to correct them. 
Anyway, any sequence of $4$-fold 
canonical flips terminates:~Theorem \ref{can2}. 
\end{abstract}

\maketitle

\section{Introduction}
Professor 
Alexeev pointed out that Lemma 2.1 in \cite{termination}, 
which is a copy of \cite[(4.12.2.1)]{FA}, 
is wrong. 
Therefore, the weighted version of difficulty 
$d_{S,b}(X,B)$ in Definition 2.3 in \cite{termination} 
is infinite if $b< \max_j\{{b_j}\}$. So, the proof in 
\cite{termination} is nonsense. 
In this paper, we change the definition of $d_{S,b}(X,B)$ 
to make it finite when $(X,B)$ is canonical and 
$B$ has no reduced components, that is, 
the round down $\llcorner B\lrcorner=0$. 
Roughly speaking, in \cite[Definition 2.3]{termination} 
we exclude valuations 
obtained by one blow-up along generic points of 
codimension two subvarieties 
when we count valuations with small discrepancies. 
In this paper, we exclude valuations whose centers 
are codimension two subvarieties with good properties.  
By this change, the new version of $d_{S,b}(X,B)$ defined 
in Definition \ref{diff} becomes finite and 
the arguments in \cite{termination} work without any changes. 
Proposition \ref{sig} is a key result in this paper. 
Note that the problems in \cite{termination} 
are not in the arguments but in the definitions. 
As mentioned above, we have to assume 
$\llcorner B\lrcorner =0$ to make $d_{S,b}(X,B)$ finite. 
Thus the main theorem:~Theorem 1.1 in \cite{termination} 
becomes slightly weaker. 
However, this assumption is harmless for 
applications if we use the special termination theorem 
(see \cite{st}). 
For the precise statements of 
the termination theorems, 
see Theorems \ref{can} and \ref{can2} below. 
Anyway, {\em{any sequence of $4$-fold canonical flips 
terminates}}. 

\begin{ack} 
I would like to thank Professor 
Valery Alexeev, 
who 
pointed out an error in \cite{termination}, gave me 
comments, and obtained 
the same correction independently. 
This paper was written in the Institute for 
Advanced Study. 
I am grateful to it for its hospitality. 
I was partially supported by a grant from the National Science 
Foundation:~DMS-0111298. 
After I wrote the preliminary version of this paper, 
I received many useful comments from Professor 
Kenji Matsuki. I would like to thank him. 
\end{ack}

We summarize the contents of 
this paper. 
In Section \ref{sec2}, we describe the errors in 
discrepancy lemmas in \cite[4.12]{FA}. 
In Section \ref{new}, we formulate a new 
discrepancy lemma. Proposition \ref{sig} 
is the main result in this paper. 
In Section \ref{sec3}, we 
explain how to modify the definition of 
the weighted version of difficulty. 
Section \ref{saigo} is devoted to the statements of 
the termination of $4$-fold canonical flips. 
We will use the same notation as in \cite{termination} throughout 
this paper. 

\section{Errors in discrepancy lemmas}\label{sec2}

The following example contradicts 
\cite[Lemma 2.1]{termination}, which is a copy of 
\cite[(4.12.2.1)]{FA}. 

\begin{exa}\label{ex1}
Let $X=\mathbb P^2$, $B=\frac{2}{3}L$, where 
$L$ is a line on $X$. Let $P$ be any point on $L$. 
First, blow up $X$ at $P$. 
Then we obtain an exceptional divisor $E_P$ such 
that $a(E_P,X,B)=\frac{1}{3}$. 
Let $L'$ be the strict transform of $L$. 
Next, take a blow-up at $L'\cap E_P$. 
Then we obtain an exceptional divisor $F_P$ whose 
discrepancy $a(F_P,X,B)=\frac{2}{3}$. 
Note that this $F_P$ is {\em{essential}} in the notation 
in \cite[Definition 2.1]{termination}. 
On the other hand, it is easy to see that 
$\xdiscrep(X,B)=\frac{1}{3}$. Thus, 
$\min\{1, 1+\xdiscrep(X,B)\}=1$. 
\end{exa}

The proof of \cite[(4.12.2.1)]{FA} depends on \cite[(4.12.1.3)]
{FA}, which is obviously wrong by Example \ref{ex1} 
above. We need some extra assumption. 
It is not difficult to see that \cite[(4.12.1.3)]{FA} is 
true if we assume that $b_j\leq \frac{1}{2}$ for all $j$. 
We write the precise statement for the reader's convenience. 
This is essentially the same as 
\cite[Corollary 3.2 (iii)]{flops} (see Remark \ref{haa} 
below).  

\begin{lmm}\label{3} 
Let $Y$ be a smooth variety with a $($not necessarily 
effective$)$ $\mathbb Q$-divisor $B=\sum_i b_i B_i$ such 
that $\sum _i B_i$ has simple normal crossings, 
$B_i$ is a prime divisor for every $i$, $B_k\ne B_l$ for 
$k\ne l$, and that 
$b_i\leq \frac{1}{2}$ for all $i$. 
Assume that $b_k+b_l\leq 0$ whenever $B_k$ and $B_l$ intersects. 
If $\nu$ is an algebraic valuation with small center on $Y$ such 
that $a(\nu, Y, B)<1$ then $\nu$ is obtained by blowing up 
the generic point of a subvariety $W\subset Y$ such 
that $\codim _Y W=2$, only one of the 
$B_k$ {\em{(}}say $B_{k_0}${\em{)}} 
contains $W$ and $b_{k_0}>0$. 
\end{lmm} 

\begin{rem}
In Example \ref{ex1}, we put $D=dL$. 
Then $a(E_P, X,D)=1-d$. Thus, the coefficient 
of $E_P$ (resp.~$D'$, the strict transform of $D$) is $d-1$ 
(resp.~$d$). 
Thus, $(d-1)+d\leq 0$ if and only if 
$d\leq \frac{1}{2}$. 
This computation shows that we have to assume 
$b_j\leq \frac{1}{2}$ for all $j$ in Lemma \ref{3}. 
\end{rem}

Thus we obtain the following lemma, which is 
a correction of \cite[(4.12.2.1)]{FA}. The proof is an exercise. 
Note that 
\cite[(4.12.2.2)]{FA} is contained in 
\cite[Proposition 2.36 (2)]{km}. 
We do not need $d_j\leq \frac{1}{2}$ for \cite[(4.12.2.2)]{FA}. 

\begin{lmm}\label{3-1}
Let $X$ be a normal variety and $D=\sum _jd_j D_j$ an 
effective $\mathbb Q$-divisor on $X$ such 
that $K_X+D$ is $\mathbb Q$-Cartier, 
where $D_j$ is a prime divisor for every $j$ and 
$D_k\ne D_l$ for $k\ne l$. 
Assume that $\xdiscrep (X,D)\geq -\frac{1}{2}$ and 
$d_j\leq \frac{1}{2}$ for all $j$. 
Let $\nu$ be an algebraic valuation with small 
center on $X$. 
Then there is a finite set of valuations $\{\nu_i\}$ 
such that if 
$$
a(\nu, X, D)<\min \{ 1, 1+\xdiscrep (X,D)\} \ \text{and} \ \ 
\nu\notin \{\nu_i\}
$$ 
then $\nu$ is obtained from blowing up the generic point 
of a subvariety 
$W\subset D\subset X$ such that 
$D$ and $X$ are generically smooth along $W$ $($and thus 
only one of the $D_j$ contains $W$$)$ and 
$\dim W=\dim X-2$. 
\end{lmm}

Unfortunately, 
Lemmas \ref{3} and \ref{3-1} are useless for our purpose. 
The assumption that $b_j\leq \frac{1}{2}$ for 
all $j$ is too strong. Proposition \ref{sig} below 
seems to be 
a better formulation. 

\begin{rem}\label{haa} 
Note that there are no problems in \cite[Corollary 3.2 (iii)]{flops} 
since Koll\'ar assumed $c>-\frac {1}{2}$ 
(for the 
notation, see Corollary 3.2 in \cite{flops}). 
The assumption $c>-\frac{1}{2}$ is in \cite[Corollary 3.2 (ii)]
{flops}. 
Lemma 2.2 in \cite{matsu} is almost an exact copy of Corollary 
3.2 in \cite{flops}. 
Therefore, \cite[Lemma 2.2]{matsu} is also correct. 
Matsuki gave me a comment about the remark 
which he made in \cite[Lemma 2.2 (ii)]{matsu} and 
which is not in \cite[Corollary 3.2]{flops} 
"(actually $>-1$ is enough for the conclusion)". 
This has to be understood that if we have the assumption $0\geq 
c>-1$, then the conclusion for (ii) holds (for the 
proof, see \cite[Proposition 2.36 (2)]{km}), and 
NOT that the conclusion of (iii) holds (as Example 
\ref{ex1} above is an obvious counter-example then). 
Thus, with the understanding that the assumptions are 
accumulative and not independent, 
it seems that the statements of the Corollary 3.2 in \cite{flops} 
and Lemma 2.2 in \cite{matsu} are correct and 
that the proof does not need any modifications. 
Therefore, the problems are not in 
\cite{flops} nor in \cite{matsu}, but 
in \cite[(4.12.1.3)]{FA}. 
For the finiteness of $d_N(X,D)$ in \cite[4.12.3 Definition]{FA}, 
we do not need \cite[(4.12.2.1)]{FA}. 
The statement \cite[(4.12.2.2)]{FA}, which is 
true by \cite[Proposition 
2.36 (2)]{km}, is sufficient. So, the error in \cite[(4.12.1.3)]{FA} 
causes no serious 
troubles in \cite[Chapter 4]{FA}.  

%Note that \cite[Corollary 3.2 (iii)]{flops} is wrong by Example 
%\ref{ex1}. 
%So, \cite[Lemma 2.2]{matsu} and \cite[(4.12.1.3)]{FA}, 
%which are copies of \cite[Corollary 3.2 (iii)]{flops}, are 
%incorrect. We think that a right formulation is \cite
%[Proposition 2.36 (2)]{km}. We note 
%that \cite[Lemma 3.3]{flops} and 
%\cite[Lemma 2.3]{matsu} are correct 
%since $\epsilon$ is sufficiently small (see Lemmas \ref{3} 
%and \ref{3-1} below). Anyway, the termination theorems proved 
%in \cite{flops}, \cite{matsu}, and \cite[Chapter 4]{FA} 
%are true. 
\end{rem}

\section{New discrepancy lemma}\label{new} 

The following proposition is a key result in this paper. The 
proof is essentially the same as one of \cite[(4.12.2.1)]{FA}. 
We give a proof for the reader's convenience. 

\begin{prp}\label{sig}
Let $X$ be a normal variety and 
$B=\sum _i b_iB_i$ a $\mathbb Q$-divisor on $X$ with 
$\llcorner B\lrcorner\leq 0$, where $B_i$ is a prime 
divisor for every $i$ and $B_k\ne B_l$ for $k\ne l$. 
Assume that $K_X+B$ is $\mathbb Q$-Cartier and 
$\xdiscrep (X,B)>-1$. Note 
that $(X,B)$ is called a sub klt pair in some literatures. 
Let $\nu$ be an algebraic valuation with 
small center on $X$. 
Then there is a finite set of valuations $\{\nu_i\}$ such 
that if 
$$
a(\nu, X, D)<\min \{ 1, 1+\xdiscrep (X,D)\} \ \text{and} \ \ 
\nu\notin \{\nu_i\}
$$ 
then $V:=\Center _X\nu \subset B \subset X$, 
$B$ and $X$ are generically smooth along $V$, 
$\dim V=\dim X-2$, and only one of the $B_k$ 
{\em{(}}say $B_{k_0}${\em{)}} 
contains $V$ and $b_{k_0}>0$. 
\end{prp}
\begin{proof}
First, we take a log resolution $f:Y\longrightarrow X$ 
as in \cite[Proposition 2.36]{km}. 
Thus, we have $f^*(K_X+B)=K_Y+A-C$, 
where $A$ and $C$ are both effective 
divisors with the following properties: 
\begin{itemize}
\item[(i)] $A=\sum_{a_i>0}a_iA_i$ and 
$C=\sum_{c_j\geq 0}c_jC_j$ have no common irreducible 
components, 
\item[(ii)] $\Exc (f)\cup \xSupp f^{-1}_*B=\sum _i A_i\cup 
\sum _j C_j$, and 
\item[(iii)] $\sum _i A_i\cup \sum _j C_j$ is a 
simple normal crossing divisor and $\sum_i A_i$ is smooth.  
\end{itemize}
Note that $c_j$ may be zero and that 
$A=f^{-1}_{*}B+D$, where $D$ is an effective 
$\mathbb Q$-divisor such that $\xSupp D\,\cap \,\xSupp 
f^{-1}_{*}B=\emptyset$. 
Next, if $E$ is an exceptional divisor 
over $Y$ such that $a(E, Y, A-C)<1$, 
then $V:=\Center _Y E \subset A\subset Y$ and 
$\dim V=\dim Y-2$ by the following 
lemma:~Lemma \ref{45}. 
We note that in general $a(E, Y, F')\leq a(E, Y, F)$ if 
$F'\geq F$ for any valuation $E$. 
If $V$ is contained in $D$, then 
$a(E,Y, A-C)\geq 1+\xdiscrep (X,B)$. 
Finally, the number of the exceptional 
divisors over $Y$ whose centers are in $f^{-1}_{*}B\cap 
C$ with $a(\cdot, Y, A-C)<1$ is finite (see Lemma \ref{45} 
below), and it is obvious that 
the number of $f$-exceptional divisors is finite. 
Thus, we obtain the required finite set of valuations 
$\{\nu_i\}$. 
\end{proof}

\begin{lmm}\label{45}
Let $Y$ be a smooth variety and $H=dP$, where 
$P$ is a smooth prime divisor on $Y$ and $0<d<1$. 
Then $\xdiscrep (Y,H)=1-d$. 
If $a(E,Y,H)<1$ for an exceptional divisor $E$ over $Y$, 
then $\Center _Y E$ is 
a codimension two subvariety of $Y$ such that 
$\Center _Y E \subset P\subset Y$. 

Let $W$ be a codimension two subvariety of $Y$ such that 
$W\subset P\subset Y$. 
Then there are only finitely many algebraic valuations 
$\nu$'s with the following properties{\em{:}} 
\begin{itemize}
\item[$(1)$] $a(\nu, Y, H)<1$, 
\item[$(2)$] $\Center _Y \nu =W$. 
\end{itemize}
Furthermore, 
$\nu$ attains the minimum, that is, $a(\nu, Y, H)=1-d$, 
if and only if $\nu$ is obtained by blowing up $Y$ along 
$W$. 
\end{lmm}
\begin{proof} 
This follows from easy computations. See 
\cite[Lemmas 2.45 and 2.29]{km}. 
\end{proof}

\section{How to define a weighted difficulty}\label{sec3} 

We introduce the notion of {\em{significant}} divisors. 
Proposition \ref{sig} and Lemma \ref{3} 
imply that the notion of significant divisors are 
much better than one of {\em{essential}} divisors 
in \cite[Definition 2.3]{termination} for our purpose. 

\begin{dfn} 
Let $(X, B)$ be a canonical pair. 
We say that an exceptional divisor $E$ 
(over $X$) is {\em{significant}} 
unless $W=\Center _X E$ is a subvariety 
$W\subset B\subset X$ 
such that $B$ and $X$ are generically smooth along $W$ (and thus only one 
of the irreducible components of $\xSupp B$ 
contains $W$) and 
$\dim W=\dim X-2$. 
\end{dfn}

The following corollary is obvious by Proposition \ref{sig}. 
We will use this to define a weighted version of difficulty. 

\begin{crl}\label{coro}
Let $(X,B)$ be a canonical pair with 
$\llcorner B\lrcorner =0$. 
Then we have 
$$
\sharp\{\,E \,|\ E \ {\text{is significant and }} a(E, X, B)<1\}<\infty. 
$$
\end{crl}

\begin{rem} Let $(X,B)$ be a canonical pair. 
Assume that $\llcorner B\lrcorner\ne 0$. 
Let $f:Y\longrightarrow X$ be a log resolution of $(X,B)$ 
with $f^*(K_X+B)=K_Y+\sum_i a_iE_i$ such 
that $\sum_i E_i=\Exc(f)\cup \xSupp f^{-1}_{*}B$. 
We can assume that 
$a_0=1$. 
If $E_0$ intersects $E_1$ such that  $0\leq 
a(E_1, X,B)=-a_1<1$ and $\codim _X f(E_0\cap E_1)\geq 3$, 
then we 
have infinitely many significant divisors whose centers are 
$f(E_0\cap E_1)$ with 
$a(\cdot, X, B)=-a_1$ by suitable blowing-ups 
whose centers are over $E_0\cap E_1$. 
\end{rem}

We define a weighted version of difficulty. 
To define this, we have to assume that the 
boundary divisor has no reduced components. 

\begin{dfn}[(A weighted version of difficulty)]\label{diff}
Let $(X,B)$ be a pair with only canonical singularities, 
where $B=\sum _{j=1}^{l}b_jB^j$ 
with $0<b_1<\cdots<b_l< 1$ 
and $B^j$ is a reduced divisor for every $j$. 
We note that $B^j$ is not necessarily 
irreducible and that we assume $b_l<1$. 
If $(X,B)$ has only terminal singularities, then 
$\llcorner B\lrcorner=0$. Thus the assumption 
$b_l<1$ always holds for terminal pairs. 
We put $b_0=0$, and $S:=\sum _{j\geq 0}b_j\mathbb Z_{\geq 0}
\subset \mathbb Q$. 
Note that $S=0$ if $B=0$. 
We set 
$$
d_{S,b}(X,B):=\sum _{\xi\in S, \xi\geq b}\sharp 
\{ E | E \ \text{is significant and} \ a(E,X,B)<1-\xi \}. 
$$ 
Then $d_{S,b_j}(X,B)$ is finite by Corollary \ref{coro}.  
\end{dfn}

\section{Statements of the termination theorems}\label{saigo}

Now the proof in \cite[\S 3]{termination} works 
without any changes only 
if we replace the word "essential" with 
"significant". Thus we obtain the following theorem, 
which is slightly weaker than the original 
theorem:~Theorem 1.1 in \cite{termination}. 

\begin{thm}\label{can} 
Let $X$ be a normal projective $4$-fold and 
$B$ an effective $\mathbb Q$-divisor such 
that $(X,B)$ is canonical and $\llcorner B \lrcorner=0$. 
Consider a sequence of log flips 
starting from $(X,B)=(X_0,B_0)${\em{:}}
$$
\begin{matrix} 
(X_0,B_0) & \dashrightarrow & {(X_1,B_1)} & \dashrightarrow & {(X_2,B_2)} 
&\dashrightarrow\cdots\\
{\ \ \ \ \searrow} & \ &  {\swarrow}\ \  {\searrow} & \ &  
{\swarrow\ \ } &\\
 \ & Z_0 & \  \ & Z_1 & \ & &, 
\end{matrix}
$$
where $\phi_i:X_i\longrightarrow Z_i$ is a contraction and 
${\phi_i}^{+}:{X_i}^{+}=X_{i+1}\longrightarrow Z_i$ is the log flip.  
Then this sequence terminates after finitely many steps. 
\end{thm} 

As we pointed out before, $\llcorner B\lrcorner=0$ if 
$(X,B)$ has only terminal singularities. 
Under the assumption that 
the varieties are $\mathbb Q$-factorial 
and all the flipping contractions have the relative Picard number 
one, we obtain the following theorem by 
using the special termination theorem. 
These assumptions are harmless for applications. 

\begin{thm}\label{can2}
Let $X$ be a normal projective $4$-fold and 
$B$ an effective $\mathbb Q$-divisor such 
that $(X,B)$ is canonical. 
Assume that $X$ is $\mathbb Q$-factorial. 
Consider a sequence of log flips 
starting from $(X,B)=(X_0,B_0)${\em{:}}
$$
\begin{matrix} 
(X_0,B_0) & \dashrightarrow & {(X_1,B_1)} & \dashrightarrow & {(X_2,B_2)} 
&\dashrightarrow\cdots\\
{\ \ \ \ \searrow} & \ &  {\swarrow}\ \  {\searrow} & \ &  
{\swarrow\ \ } &\\
 \ & Z_0 & \  \ & Z_1 & \ & &, 
\end{matrix}
$$
where $\phi_i:X_i\longrightarrow Z_i$ is a contraction and 
${\phi_i}^{+}:{X_i}^{+}=X_{i+1}\longrightarrow Z_i$ is the log flip.  
We further assume that the relative Picard 
number $\rho(X_i/Z_i)=1$ for every $i$. 
Then this sequence terminates after finitely many steps. 
\end{thm}
\begin{proof}
By applying the special termination theorem (see \cite{st}) 
and 
shifting the index, we 
can assume that the flipping and flipped loci are disjoint from 
$\llcorner B_i\lrcorner$ for every $i$. 
So, we can replace $B_i$ with its fractional part. 
Thus this sequence terminates by 
Theorem \ref{can}. 
\end{proof}

\begin{rem}
The final remark in \cite{termination} should be removed. 
In \cite{ssfli}, we only need the termination of 
$4$-fold semi-stable terminal 
flips. See Definition 2.3 in \cite{ssfli}. 
Therefore, Theorems \ref{can} and \ref{can2} are 
sufficient for \cite{ssfli}. 
\end{rem}
%%%%%%%%%%%%%%%%%%%%%%%%%%%%%%%%%


\begin{thebibliography}{KM}
%
% The \bibitem commands: 
% Please follow "Notice to Authors" for referencing.  You 
% must specify bold and italic fonts yourself. 
%

\bibitem[F1]{termination} 
O.~Fujino, Termination of 4-fold canonical flips, 
Publ. Res. Inst. Math. Sci. {\textbf{40}} (2004), no.1, 
231--237. 

\bibitem[F2]{ssfli} 
O.~Fujino, 
On termination of 4-fold semi-stable log flips, 
to appear in Publ. Res. Inst. Math. Sci.  

\bibitem[F3]{st} 
O.~Fujino, Special termination and reduction theorem, 
to be contained in the book prepared by 
A.~Corti et al. 

\bibitem[K]{flops} 
J.~Koll\'ar, Flops, Nagoya Math. J. {\textbf{113}} 
(1989), 15--36.

\bibitem[KM]{km} 
J.~Koll\'ar, and S.~Mori, 
{\em{Birational geometry of algebraic varieties,}} 
Cambridge Tracts in Mathematics, {\textbf{134}}. 
Cambridge University Press, Cambridge, 1998. 

\bibitem[K$^+$]{FA} 
J.~Koll\'ar et al., {\em{Flips and Abundance for algebraic threefolds}}, 
Ast\'erisque \textbf{211}, (1992).  

\bibitem[M]{matsu} 
K.~Matsuki, 
Termination of flops for $4$-folds, 
Amer. J. Math. {\textbf{113}} (1991), no. 5, 835--859.

\end{thebibliography}
\end{document}